\documentclass[12pt]{amsart}
\usepackage{amsmath,amssymb,amsfonts,amsthm,amsopn,dsfont}
\usepackage{graphics}

\newcommand{\ft}{Fourier transform}

\newcommand{\fpc}{\Fur L^p(\R^d)_{\rm comp}}

\newtheorem{tm}{Theorem}[section]
\newtheorem{lemma}[tm]{Lemma}

\newtheorem{theorem}{Theorem}[section]
\newtheorem{corollary}[theorem]{Corollary}
\newtheorem{definition}[theorem]{Definition}

\newtheorem{proposition}[theorem]{Proposition}

\newcommand{\beqa}{\begin{eqnarray*}}
\newcommand{\eeqa}{\end{eqnarray*}}

\DeclareMathOperator*{\supp}{supp}

\newcommand{\field}[1]{\mathbb{#1}}
\newcommand{\bR}{\field{R}}        
\newcommand{\bZ}{\field{Z}}        
\def\la{\lambda}

 \def\cF{\mathcal{F}}              
 \def\cS{\mathcal{S}}

 \def\cC{\mathcal{C}}

\def\a{\aleph}

\def\rd{\bR^d}

\def\rdd{{\bR^{2d}}}

\def\zd{\bZ^d}

\def\intrd{\int_{\rd}}

\def\R{\right)}

\def\<{\left<}
\def\>{\right>}

\def\mv1{M_v^1}

\def\phas{(x,\o )}
\def\mn{(m,n)}
\def\mn'{(m',n')}

\def\o{\eta}
\def\a{\alpha}
\def\b{\beta}

\def\N{\mathbb{N}}
\def\R{\mathbb{R}}
\def\Ren{\mathbb{R}^d}

\def\Fur{\mathcal{F}}

\def\f{\varphi}

\def\Sn2{S_{2}(L^{2}(\Ren))}
\def\S1{S_{1}(L^{2}(\Ren))}
\def\sig00{\sigma_{0,0}}

\def\la{\langle}
\def\ra{\rangle}

\begin{document}

\begin{abstract}
We carry on the study of Fourier
integral operators of
H{\"o}r\-man\-der's type acting on the
spaces $\fpc$, $1\leq p\leq\infty$, of
compactly supported distributions whose
Fourier transform is in $L^p$. We show
that the sharp loss of derivatives for
such an operator to be bounded on these
spaces is related to the rank $r$ of
the Hessian of the phase $\Phi(x,\eta)$
with respect to the space variables
$x$. Indeed, we show that operators of
order $m=-r|1/2-1/p|$ are bounded on
$\fpc$, if the mapping
$x\longmapsto\nabla_x\Phi(x,\eta)$ is
constant on the fibers, of codimension
$r$, of an affine fibration.
\end{abstract}

\title[Fourier integral
Operators on $\Fur
L^p$]{Boundedness of Fourier
integral operators on Fourier
Lebesgue spaces and affine
fibrations}\author{Fabio
Nicola}
\address{Dipartimento di Matematica,
Politecnico di Torino, corso Duca degli
Abruzzi 24, 10129 Torino, Italy}
\email{fabio.nicola@polito.it}
\thanks{}
\maketitle

\section{Introduction}
Consider the spaces $\fpc$ of
compactly supported
distributions whose Fourier
transform is in $L^p(\R^d)$,
with the norm $\|f\|_{\Fur
L^p}=\|\hat{f}\|_{L^p}$. In
\cite{cordero-nicola-rodino}
we studied the boundednss on
these spaces of H\"ormander's
type Fourier integral
operators (FIOs), namely
operators of the form
\begin{equation}\label{FIO}
Tf(x)=\int e^{2\pi
i\Phi(x,\eta)}
\sigma(x,\eta)\hat{f}(\eta)\,d\eta.
\end{equation}
Here the \ft\, of
$f\in\cS(\rd)$ is normalized
to be ${\hat
  {f}}(\o)=\int
f(t)e^{-2\pi i t\o}dt$.
 The
symbol $\sigma$ is in
$S^m_{1,0}$, the
H{\"o}rmander's class of
order $m$. Namely, $\sigma\in
\mathcal{C}^\infty(\R^{2d})$
and satisfies
\begin{equation}\label{symb}
|\partial^\alpha_x\partial^\beta_\eta
\sigma(x,\eta)|\leq
C_{\alpha,\beta}
\langle\o\rangle^{m-|\beta|},
\quad \forall
(x,\o)\in\R^{2d},
\end{equation}
where, as usual,
$\langle\eta\rangle:=(1+|\eta|^2)^{1/2}$.
We also suppose that $\sigma$
has compact support with
respect to $x$.\par The phase
$\Phi(x,\eta)$ is
real-valued,
 positively
 homogeneous of degree 1 in $\eta$, and smooth on
 $\R^d\times(\R^d\setminus\{0\})$.  We assume
 $\Phi(x,\eta)$ defined on
 an open
subset
 $\Lambda\subset \R^d\times(\R^d\setminus\{0\})$,
 conic in dual variables,
containing the closure of the set
\[
\Lambda'=\{(x,\eta)\in\R^d\times(\R^d\setminus\{0\}):
(x,\lambda\eta)\in{\rm
supp}\,\sigma\ \textrm{for
some}\ \lambda>0\},
\]
 in
$\R^d\times(\R^d\setminus\{0\})$.
We also assume the
non-degeneracy condition
\begin{equation}\label{nondeg}
{\rm det}\,
\left(\frac{\partial^2\Phi}{\partial
x_i\partial \eta_l}\Big|_{
(x,\eta)}\right)\not=0\quad \forall
(x,\eta)\in \Lambda.
\end{equation}
It is easy to see that such
an operator maps the space
$\mathcal{S}(\R^d)$ of
Schwartz functions into the
space
$\mathcal{C}^\infty_0(\R^d)$
of test functions
continuously. See
\cite{hormander,treves} for
the general theory of FIOs,
and \cite{sogge93,stein93}
for results in $L^p$.

\par The main result of \cite{cordero-nicola-rodino} states
that, if
\[
m\leq-d\left|\frac{1}{2}-\frac{1}{p}\right|,
\]
 the operator $T$
above, initially defined on
$\mathcal{C}^\infty_0(\R^d)$, extends
to a bounded operator on $\fpc$,
whenever $1\leq p<\infty$. For
$p=\infty$, $T$ extends to a bounded
operator on the closure of
$\mathcal{C}^\infty_0(\R^d)$ in $\Fur
L^\infty(\R^d)_{\rm comp}$. Moreover
this loss of derivatives was proved to
be generally sharp in any dimension,
even for phases linear in $\eta$. This
result can be regarded as an
investigation of the minimal loss of
derivatives occurring in the
Beurling-Helson theorem
\cite{beurling53,lebedev94} (see also
\cite{kasso07,rstt}).\par However, it
was shown in \cite{rstt} that, in the
case of phases linear in $x$, local
$\Fur L^p$-boundedness holds without
loss of derivatives, i.e. for $m=0$.
This suggests the possibility of
intermediate thresholds, depending on
the rank of Hessian
$d^2_x\Phi(x,\eta)$.\par A similar
phenomenon is well-known for
$L^p$-boundedness. Namely, a celebrated
result by Seeger, Sogge and Stein
\cite{seeger-sogge-stein} shows that
$T$ is bounded on $L^p$, $1<p<\infty$,
if $m\leq -(d-1)|1/2-1/p|$ (see also
\cite{stein93,tao}). Moreover, if the
rank of the Hessian
$d^2_\eta\Phi(x,\cdot)$ is $\leq r$,
then the threshold goes up to
$-r|1/2-1/p|$, provided a certain
smooth factorization condition is
satisfied. Although that condition is
not necessary for boundedness to hold,
it turns out to be essential in the
proof, given in
\cite{seeger-sogge-stein}, when the
rank of $d^2_\eta\Phi(x,\eta)$ is
allowed to drop, and its relaxation is
an open problem. The main reference on
this topic is Ruzhansky's survey
\cite{ruz2} and book
\cite{ruzhansky-book}, where the case
of complex-valued phases is also
considered. See also \cite{ruz1}.\par
In this paper we present a variant of
the smooth factorization condition
which is relevant when dealing with
$\Fur L^p$ spaces. Then we show, under
that condition, that in fact the above
threshold for local $\Fur
L^p$-boundedness can go up.\par
\begin{definition}[Spatial smooth factorization
condition]\label{ssf} Let
$0\leq r\leq d$ and suppose
that for every
$(x_0,\eta_0)\in\Lambda$,
$|\eta_0|=1$, there exists an
open neighborhood $\Omega$ of
$x_0$ and an open
neighborhood
$\Gamma'\subset\mathbb{S}^{d-1}$
of $\eta_0$, with
$\Omega\times\Gamma'\subset\Lambda$,
satisfying the following
condition. For every
$\eta\in\Gamma'$ there exists
a smooth
 fibration of $\Omega$, smoothly depending on
 $\eta$ and with affine fibers of
 codimension $r$,
 such that
 $\nabla_x\Phi(\cdot,\eta)$
 is constant on every
 fiber.\footnote{To be precise,
 by ``a fibration of $\Omega$,
 smoothly depending on $\eta\in\Gamma'$ and with fibers
 of codimension $r$" we mean that
  a smooth function
 $\Pi:\Omega\times\Gamma'\to\R^d$ is given,
  with $d_x\Pi$ having constant rank $r$. The fibers are the level sets of the mapping $\Pi(\cdot,\eta)$.}
 \end{definition}
Observe that this condition
implies the Hessian
$d^2_x\Phi(x,\eta)$ to have
rank $\leq r$. Moreover it is
always satisfied if $r=d$ or
if $d^2_x\Phi(x,\eta)$ has
constant rank $r$ (in
particular, for phases linear
in $x$, corresponding to
$r=0$).

\begin{theorem}\label{maintheorem2}
Assume the above hypotheses
on the symbol $\sigma$ and
the phase $\Phi$. Assume,
moreover, that $\Phi$
satisfies the spatial smooth
factorization condition
(Definition \ref{ssf}) for
some $r$. If
\begin{equation}\label{soglia}
m\leq-r\left|\frac{1}{2}-\frac{1}{p}\right|,
\end{equation}
then the corresponding FIO T,
initially defined on
$\mathcal{C}^\infty_0(\R^d)$,
extends to a bounded operator
on $\fpc$, whenever $1\leq
p<\infty$. For $p=\infty$,
$T$ extends to a bounded
operator on the closure of
$\mathcal{C}^\infty_0(\R^d)$
in $\Fur L^\infty(\R^d)_{\rm
comp}$.
\end{theorem}
It is possible to see that
the threshold in
\eqref{soglia} is
 sharp in any dimension
$d\geq1$, even for phases
$\Phi(x,\eta)$ which are
linear in $\eta$. Namely,
 in dimension $d$, consider the phase
$\Phi(x,\eta)=\sum_{k=1}^r\varphi(x_k)\eta_k+\sum_{k=r+1}^d
x_k\eta_k$, where
$\varphi:\R\to\R$ is a
diffeomorphism,
 with $\f(t)=t$ for $|t|\geq 1$ and whose
 restriction to $(-1,1)$ is non-linear.
 Consider then the symbol
 \[
\sigma\phas=G(x)\la
\eta\ra^{{m}},\quad
\mbox{with}\,\,\,G\in\cC_0^\infty(\R^d),\,\,\,
G\equiv1\,\,\mbox{on}\,\,[-1,1]^d.
\]
Let moreover $1\leq p\leq2$.
Then Theorem
\ref{maintheorem2} and an
easy variant of the arguments
in \cite[Section
6]{cordero-nicola-rodino}
show that the corresponding
operator $T:\Fur
L^p(\R^d)_{\rm comp}\to \Fur
L^p(\R^d)_{\rm comp}$ is
bounded if and only if $m$
satisfies \eqref{soglia}. By
duality one can construct
similar examples for
$2<p\leq\infty$.\par The
proof of the result in
\cite{cordero-nicola-rodino},
corresponding to Theorem
\ref{maintheorem2} with $r=d$
(hence the spatial smooth
factorization condition is
automatically satisfied) used
tools from Time-frequency
analysis, relying on our
previous work \cite{fio1}.
Instead, the proof of Theorem
\ref{maintheorem2} is
inspired by more classical
arguments in
\cite{seeger-sogge-stein}.
Indeed, we will conjugate the
operator $T$ above with the
Fourier transform, obtaining
the FIO
\[
\tilde{T} f(x)=\Fur\circ
T\circ\Fur^{-1} f(x)=\iint e^{2\pi
i(\Phi(\eta,y)-x\eta)}\sigma(\eta,y)
f(y)\, dy\,d\eta,
\]
for which $L^p$-boundedness
has to be proved. However
notice that we cannot apply
to $\tilde{T}$ the well-known
results of the classical
(i.e.\ local) $L^p$-theory.
Indeed, the corresponding
symbol is no longer compactly
supported with respect to $y$
and the phase is no longer
homogeneous with respect to
$\eta$. Among other things,
it does not satisfy the
(frequency) smooth
factorization condition of
\cite{seeger-sogge-stein}. In
fact, for operators
$\tilde{T}$ of this special
type we will prove results in
$L^p$ in the limiting cases
$p=1,\infty$ too, which are
generally false, for example,
for the operator $T$ above.
\par Our strategy consists in
splitting $\tilde{T}$ in
dyadic pieces via a
Littlewood-Paley
decomposition of the {\it
physical} domain and then
each dyadic operator is
further split in a certain
number of FIOs with symbols
localized in thin boxes of
the {\it frequency} domain,
and phases essentially linear
in $\eta$. By comparison,
notice that in the classical
$L^p$-theory one performs a
dyadic decomposition and then
a second decomposition, both
in the frequency domain
(\cite{seeger-sogge-stein,stein93});
moreover, the geometry of our
second decomposition is
different from that in
\cite{seeger-sogge-stein}.\par
 This discussion also shows
 that Theorem
 \ref{maintheorem2} can be
 read as a {\it global} boundedness
 result on $L^p$ for the
 operator $\tilde{T}$,
 and hence it partially
 intersects some recent results in
\cite{coriascoruz} (see also
some examples in
\cite{cordero-nicola-rodino2}
and, for the case $p=2$,
\cite{coriasco,ruzhsugimoto}).
In fact, sharp results are
obtained there for general
classes of operators, but no
improvements upon them seem
to be considered, under
additional conditions like a
smooth
factorization.\par\medskip\noindent{\bf
Notation} We write $A\lesssim
B$ if $A\leq C B$ for some
constant $C>0$ which may
depend on parameters, like
Lebesgue exponents or the
dimension $d$. We write
$A\asymp B$ if $A\lesssim B$
and $B\lesssim A$. Finally,
for $R>0$, $x_0\in\R^d$, we
denote by $B_R(x_0)$ the open
ball in $\R^d$ with centre
$x_0$ and radius $R$.

\section{Preliminary results on
FIOs}\label{section3} In the
next section we will make use
of the  well-known
composition formula of a
pseudodifferential operator
and a FIO. We collected here
what is needed in the
subsequent proofs.
\par
First we recall that a
regularizing operator is a
pseudodifferential operator
\[
Rf(x)=r(x,D)f=\int e^{2\pi
ix\eta}r(x,\eta)\hat{f}(\eta)d\eta,
\]
with a symbol $r$ in the
Schwartz space
$\mathcal{S}(\R^{2d})$
(equivalently, an operator
with kernel in
$\mathcal{S}(\R^{2d})$, which
maps $\mathcal{S}'(\R^{d})$
into $ \mathcal{S}(\R^{d})$).
Then, the composition formula
for a pseudodifferential
operator and a FIO is as
follows (see, e.g.,
\cite{hormander},
\cite[Theorem
4.1.1]{mascarello-rodino},
\cite[Theorem
18.2]{Shubin91},
\cite{treves}).
\begin{theorem}\label{composition}
Let the symbol $\sigma$ and
the phase $\Phi$ satisfy the
assumptions in the
Introduction. Assume, in
addition, $\sigma(x,\eta)=0$
for $|\eta|\leq1$, if
$\Phi(x,\eta)$ is not linear
in $\eta$. Let $q(x,\eta)$ be
a symbol in $S^{m'}_{1,0}$.
Then,
\[
q(x,D)T=S+R,
\]
where $S$ is a FIO with the same phase
$\Phi$ and symbols $s(x,\eta)$, of
order $m+m'$, satisfying
\[
{\rm supp}\,s\subset {\rm
supp}\,\sigma\cap\{(x,\eta)\in\Lambda:\
(x,\nabla_x\Phi(x,\eta))\in{\rm
supp}\,q\},
\]
and $R$ is a regularizing
operator with symbol
$r(x,\eta)$ satisfying
\[
\Pi_\eta ({\rm
supp}\,r)\subset \Pi_\eta(
{\rm supp}\,\sigma),
\]
where $\Pi_\eta$ is the
orthogonal projection on
$\R^d_\eta$.\par Moreover,
the symbol estimates
satisfied by $s$ and the
seminorm estimates of $r$ in
the Schwartz space are
uniform when $\sigma$ and $q$
vary in a bounded subsets of
$S^m_{1,0}$ and
$S^{m'}_{1,0}$ respectively.
\end{theorem}

\section{Proof of the main result (Theorem \ref{maintheorem2})}

By means of a smooth cut off function
near $\eta=0$ we split the symbol
$\sigma$ of $T$ in a symbol supported
where $|\eta|\leq 4$ and a symbol
supported where $|\eta|\geq2$. Now, the
first symbol yields an operator which
is bounded on all $\Fur L^p(\R^d)_{\rm
comp}$, $1\leq p<\infty$, as well as on
the closure of $\cC^\infty_0(\R^d)$ in
$\Fur L^\infty(\R^d)_{\rm comp}$. This
was shown in \cite[Proposition
4.1]{cordero-nicola-rodino}, regardless
of the order of the operator, and
without assuming the condition
\eqref{nondeg}, nor the spatial smooth
factorization condition\footnote{There
the desired boundedness was proved on
the so-called modulation spaces $M^p$.
However it was observed that the
corresponding norm is equivalent to
that of $\Fur L^p$ for distributions
supported in a fixed compact set.}.\par
Hence we will prove the estimate
\begin{equation}\label{stimas}
\|T f\|_{\Fur L^p}\lesssim
\|f\|_{\Fur L^p},\quad
\forall
f\in\cC^\infty_0(\R^d),
\end{equation}
$1\leq p\leq\infty$, for an
operator satisfying the
assumptions of Theorem
\ref{maintheorem2} and whose
symbol $\sigma$ satisfies, in
addition,
\[
\sigma(x,\eta)=0\quad {\rm
for}\ |\eta|\leq 2.
\]
This clearly implies the
conclusion of Theorem
\ref{maintheorem2}.\par
 We
first perform a further
reduction.
 For every
$(x_0,\eta_0)\in \Lambda'$,
$|\eta_0|=1$, there exist an
open neighbourhood
$\Omega\subset\R^d$ of $x_0$,
an open conic neighbourhood
$\Gamma\subset\R^d\setminus\{0\}$
of $\eta_0$ and $\delta>0$
such that
\begin{equation}\label{i1}
|\det
\partial^2_{x,\eta}\Phi(x,\eta)|\geq\delta>0,\quad
\forall
(x,\eta)\in\Omega\times\Gamma,
\end{equation}
and
\begin{equation}\label{i2}
\forall x\in\Omega,\
\textrm{the map}\
\Gamma\ni\eta\mapsto
\nabla_x\Phi(x,\eta)\
\textrm{is a diffeomorphism
onto the range}.
\end{equation}
Hence, by a compactness
argument and a finite
partition of unity we can
assume that $\sigma$ itself
is supported in a set of the
type $\Omega'\times\Gamma$,
for some open
$\Omega'\subset\subset\Omega\subset\R^{d}$,
and
 conic open $\Gamma\subset\R^{d}\setminus\{0\}$,
  with $\Phi$ satisfying the
above conditions on
$\Omega\times\Gamma$, as well
as the spatial smooth
factorization condition
(Definition \ref{ssf}) for
$x\in\Omega$,
$\eta\in\Gamma':=\Gamma\cap\mathbb{S}^{d-1}$.\par
Now, we will prove
\eqref{stimas} with
$p=1,\infty$, for an operator
$T$ of order $m=-r/2$. Then
the desired result when
$1<p<\infty$, for operators
of order $m=-r|1/2-1/p|$,
will follow by complex
interpolation with the
well-known case
 $L^2$, (see e.g.
 \cite[page 402]{stein93})\par
 In detail, for the interpolation step we argue as follows.
 For $s\in\R$,
 denote by $\Fur L^p_s$ the
 space of temperate distributions $f$ such that
 \[
 \|f\|_{\Fur
 L^p_s}:=\left(\int
 \langle\eta\rangle^{ps}|\hat{f}(\eta)|^p\,d\eta\right)^{1/p}<\infty,
 \]
 with the obvious changes if
 $p=\infty$.
For
 every $s\in\R$, the operator
 $\langle D\rangle^s$ defines
 an isomorphism of $\Fur L^p_s$
 onto $\Fur L^p$.
   Hence, the
 operator $T=T\langle D\rangle^{-s}\langle D\rangle^s$
 is bounded $\Fur L^p_s\to \Fur L^p$
 if
 $T\langle D\rangle^{-s}$ is
 bounded on $\Fur L^p$. Observe
 moreover that $T\langle
 D\rangle^{-s}$ is a FIO
 with the same phase as $T$, and symbol
 $\sigma(x,\eta)\langle\eta\rangle^{-s}$,
 which has order $m-s$.\\
 Suppose now that the
 desired result is already
 obtained for $p=1,2$.
 Take $1<p<2$ and
consider a FIO $T$ of order
 $m=-r(1/p-1/2)$. Then, taking the above remarks into
 account, $T$ extends to a
 bounded operator
 $\Fur L^1_{m+r/2}\to \Fur L^1$ and
 $L^2_m\to L^2$. Hence, the
  boundedness on $\Fur L^p$ follows
 by complex interpolation, because, if
 $\theta\in(0,1)$ satisfies
 $(1-\theta)/1+\theta/2=1/p$,
 one has
 $(m+r/2)(1-\theta)+m\theta=0$. The
 proof for $2<p<\infty$ is
 similar.\par Of course, when
 in \eqref{soglia}
 there is a strict inequality, the desired result  follows from the
  equality-case, for an operator with order $m'<m$ has also order $m$. \par
Hence, from now on, we assume $m=-r/2$
and prove \eqref{stimas} for
$p=1,\infty$.\par The first step
consists in conjugating $T$ with the
Fourier transform. The desired results
will be proved if we verify that the
operator $\tilde{T}=\Fur\circ T\circ
\Fur^{-1}$ is continuous on $L^1$ and
on the closure of $\cC^\infty_0$ in
$L^\infty$. This operator has integral
kernel
\[
K(x,y)=\int e^{2\pi
i(\Phi(\eta,y)-x\eta)}\sigma(\eta,y)\,d\eta,
\]
which is smooth everywhere
and supported in
$\R^d\times\Gamma$. Indeed,
as anticipated in the
Introduction, the problem is
the integrability at
infinity.\par Consider now
the usual Littlewood-Paley
decomposition, but on the
physical domain. Namely,
 fix a smooth function $\psi_0(y)$
  such that $\psi_0(y)=1$
  for $|y|\leq1$ and
  $\psi_0(y)=0$ for
  $|y|\geq2$. Set
  $\psi(y)=\psi_0(y)-\psi_0(2y)$,
  $\psi_j(y)=\psi(2^{-j}y)$, $j\geq1$.
  Then
  \[
  1=\sum_{j=0}^\infty\psi_j(y),\quad
\forall y\in\R^d.
 \]
  Notice
that, if $j\geq 1$, $\psi_j$ is
supported where $2^{j-1}\leq|y|\leq
2^{j+1}$. Since $\sigma(\eta,y)=0$, for
$|y|\leq 2$, we can write the kernel
above as
\[
K=\sum_{j\geq1} K_{j},
\]
where
\[
K_j(x,y)=\int e^{2\pi
i(\Phi(\eta,y)-x\eta)}\sigma_j(\eta,y)\,d\eta,
\]
and
\[
\sigma_j(\eta,y):=\sigma(\eta,y)\psi_j(y).
\]
Observe that $\eta$ lies in the open
neighbourhood $\Omega'$. After
shrinking $\Omega'$ and $\Gamma$, if
necessary, we see from the spatial
smooth factorization condition that
there exist an open neighbourhood
$U\times V$ of $(0,0)$ in
$\R^r\times\R^{d-r}$ and a smooth
change of variables $U\times
V\ni(u,v)\longmapsto\eta_y(u,v)\in\Omega_y$,
smoothly depending on the parameter
$y\in\Gamma$ and homogeneous of degree
$0$ with respect to $y$, with
$\Omega'\subset\Omega_y\subset\Omega$,
such that the function
$\eta\longmapsto\nabla_1\Phi(\eta,y)$
in constant on each of the
$(d-r)$-dimensional (pieces of) affine
planes $u=const$. Here we denoted by
$\nabla_1$ the gradient with respect to
the first $d$ variables, i.e.
$\nabla_1\Phi(\eta,y)
=\nabla_\eta\left(\Phi(\eta,y)\right)$.
Pay attention that, because of the
previous exchange of variables, we are
using Definition \ref{ssf} with the
variables $(x,\eta)$ replaced by
$(\eta,y)$.
\par
For every $j\geq1$ we choose $u^\nu_j$,
$\nu=1,2,..., N_r(j)$, such that
$|u^\nu_j-u^{\nu'}_j|\geq C_0 2^{-j/2}$
for different $\nu,\nu'$, and such that
$U$ is covered by balls with centre
$u^\nu_j$ and radius $C_1 2^{-j/2}$. It
is easy to see that
$N_r(j)=O(2^{jr/2})$. Let then
$\eta^\nu_j=\eta_y(u^{\nu}_j,0)$.
Consider moreover a smooth partition of
unity tailored to the covering above,
namely given by smooth functions
$\chi^\nu_j(u)$, $\nu=1,2,..., N_r(j)$,
supported in the above balls, and
satisfying the estimate
\begin{equation}\label{perdita}
|\partial^\alpha_u\chi^\nu_j(u)|\lesssim
2^{j|\alpha|/2}.
\end{equation}
Accordingly we decompose the kernel
$K_j=\sum_{\nu=1}^{N_r(j)} K^\nu_j$,
where
\begin{equation}\label{integrale}
K^\nu_j(x,y)=\int e^{2\pi
i(\Phi(\eta,y)-x\eta)}\chi^\nu_j(u(y,\eta))
\sigma_j(\eta,y)\,d\eta,
\end{equation}
where the function
$u=u(y,\eta)$ is obtained by
the inverse change of
variables. Consider now the
second order Taylor expansion
of $\Phi(\cdot,y)$ at
$\eta^\nu_j$:
\[
\Phi(\eta,y)=\Phi(\eta^\nu_j,y)+
\nabla_1\Phi(\eta^\nu_j,y)(\eta-\eta^\nu_j)+
R^\nu_j(\eta,y),\] where
\begin{equation}\label{restotaylor}
R^\nu_j(\eta,y)=\frac{1}{2}\int_0^1
(1-t)
\left(d^2_1\Phi\right)(\eta^\nu_j+t(\eta-\eta^\nu_j),y)
[\eta-\eta^\nu_j,
\eta-\eta^\nu_j]\,dt.
\end{equation}
Here $d^2_1\Phi$ stands for the Hessian
of $\Phi$ with respect to the first $d$
variables, regarded as a bilinear
form.\par For fixed $j$ and $\nu$,
after a rotation we perform a splitting
$(\eta',\eta'')$ of $\eta$ such that
the vectors $\eta''$ are tangent to the
pieces of plane $u=u^\nu_j$. We
consider then the operator
\[
L^\nu_j=(1-\langle2^{-j/2}\nabla_{\eta'},2^{-j/2}
\nabla_{\eta'}\rangle)(1-\langle\nabla_{\eta''},
\nabla_{\eta''}\rangle)
\]
We will prove in the Appendix
at the end of this paper that
\begin{equation}\label{daprovare}
|(L^\nu_j)^N\left(e^{2\pi
i(\Phi(\eta^\nu_j,y)-\langle\nabla_1\Phi(\eta^\nu_j,y),\eta^\nu_j\rangle+R^{\nu}_j(\eta,y))}
\chi^\nu_j(u(y,\eta))\sigma_j(\eta,y)\right)|\leq
C_N 2^{-jr/2}.
\end{equation}
 An
integration by parts in
\eqref{integrale} then yields
\[
|K^\nu_j(x,y)|\lesssim
2^{-jr}(1+2^{-j/2}|(\nabla_1\Phi(\eta^\nu_j,y)-x)'|)^{2N}
(1+|(\nabla_1\Phi(\eta^\nu_j,y)-x)''|)^{2N},
\]
where we took into account that we
integrated on a set of measure
$\lesssim 2^{-jr/2}$.\par Hereby  one
sees that, choosing $N>d/2$,
\begin{equation}\label{int1}
\int
|K^\nu_j(x,y)|\,dx\lesssim
2^{-jr/2}.
\end{equation}
To treat the integral of the
kernel with respect to $y$,
in view of \eqref{i2} we can
perform the change of
variable
$y\longmapsto\nabla_1\Phi(\eta^\nu_j,y)$,
whose inverse jacobian
determinant
 is homogeneous of
degree $0$ and uniformly
bounded because of
\eqref{i1}. We therefore
obtain
\begin{equation}\label{int2}
\int
|K^\nu_j(x,y)|\,dy=\int_\Gamma
|K^\nu_j(x,y)|\,dy\lesssim
2^{-jr/2},
\end{equation}
 It follows from \eqref{int1}, \eqref{int2}
 and the Schur
 test (see e.g.
 \cite[Theorem 6.18]{folland})
 that the operators with
 kernel $K^\nu_j$ are bounded
 on $L^1$ and $L^\infty$, with
 norm $O(2^{-j/2})$. Summing
 over $\nu=1,...,N_r(j)$ we
 see that the operators
with kernel $K_j$
 satisfy the uniform
 estimates
 \[
 \|\int K_j(\cdot,y)f(y)\,dy\|_{L^p}\lesssim
 \|f\|_{L^p},\quad
 p=1,\infty,\quad
 \forall u\in\cS(\R^d).
 \]
 Coming back to the $\Fur
 L^p$ spaces we see that the
 operators
 \[
 T_jf(x):=\int e^{2\pi
i\Phi(x,\eta)}
\sigma(x,\eta)\psi_j(\eta)\hat{f}(\eta)\,d\eta
\]
therefore fulfill the
estimates
 \begin{equation}\label{stime}
 \|T_j f\|_{\Fur L^p}\lesssim
 \|f\|_{\Fur L^p},\quad
 p=1,\infty,\quad
 \forall u\in\cS(\R^d).
 \end{equation}
 Hence, in order to obtain the desired estimate
  for the original operator $T=\sum_{j=1}^\infty
  T_j$, we are left to summing the
 estimates \eqref{stime}
 over
 $j\geq1$.\par\medskip\noindent
 {\bf Summing estimates
 \eqref{stime} for $p=1$.}
 Observe that if $\chi$ is a
smooth function supported
where $B_0^{-1}\leq|\eta|\leq
B_0$, for some $B_0>0$, then
trivially
\begin{equation}\label{sommap1}
\sum_{j=1}^\infty\|\chi(2^{-j}
D)f\|_{\Fur L^1}\lesssim
\|f\|_{\Fur L ^1},
\end{equation} for every
$f\in\cS(\rd)$, where
$\chi(2^{-j} D)f=
\cF^{-1}[\chi(2^{-j}\cdot)\hat{f}]$.

Also, notice that, due to the
frequency localization of
$T_{j}$, the estimate
\eqref{stime} for $p=1$ can
be refined as
\[
\|T_{j} f\|_{\Fur
L^1}=\|T_{j}(\chi(2^{-j}D)
f)\|_{\Fur
L^1}\lesssim\|\chi(2^{-j}D)f\|_{\Fur
L^1},
\]
where $\chi$ is a smooth
function satisfying
$\chi(\eta)=1$ for $1/2\leq
|\eta|\leq2$ and
$\chi(\eta)=0$ for
$|\eta|\leq1/4$ and
$|\eta|\geq4$ (so that
$\chi\psi=\psi$). Summing
over $j$ this last estimate
with the aid of
\eqref{sommap1} we obtain
\[
\|Tf\|_{\Fur
L^1}\lesssim\|f\|_{\Fur L^1},
\]
which is the desired
estimate.\par\medskip\noindent
{\bf Summing estimates
 \eqref{stime} for $p=\infty$.}

Here we use the following
trivial remark. For $k\geq0$,
let $f_k\in\cS(\rd)$ satisfy
${\rm supp}\,\hat{f}_0\subset
B_2(0)$ and
\[
{\rm
\supp}\,\hat{f}_k\subset\{\eta\in\R^d:\
2^{k-1}\leq|\eta|\leq
2^{k+1}\},\quad k\geq1.
\]
Then, if the sequence $f_k$ is bounded
 in $\Fur L^\infty§(\rd)$, the series
$\sum_{k=0}^\infty {f}_k$
converges in ${\Fur
L}^\infty(\R^d)$ and
\begin{equation}\label{b0}
\|\sum_{k=0}^\infty
f_k\|_{\Fur
L^\infty}\lesssim\sup_{k\geq0}\|f_k\|_{\Fur
L^\infty}.
\end{equation}
 Indeed, at each point in the frequency domain there
  are at most two
nonzero terms in the sum.
This observation yields
\begin{align}
\|Tf\|_{\Fur
L^\infty}&=\|\sum_{k\geq0}\psi_k(D)
Tf\|_{\Fur L^\infty}\nonumber\\
&\lesssim\sup_{k\geq0}\|\psi_k(D)Tf\|_{\Fur L^\infty}\nonumber\\
&=\sup_{k\geq0}\|\sum_{j=1}^\infty\psi_k(D)T_{j}
f\|_{\Fur
L^\infty}.\label{continuazione}
\end{align}
Notice that the sequence of
symbols $\sigma_j(x,\eta)$ is
bounded in $S^{-r/2}_{1,0}$,
whereas the sequence of
symbols $\psi_k(\eta)$ is
bounded in $S^0_{1,0}$.
\par Applying  Theorem
\ref{composition} to each
product $\psi_k(D)T_{j}$, we
have
\begin{equation}\label{psitsr}
\psi_k(D)T_{j}=S_{k,j}+R_{k,j},
\end{equation}
where $S_{k,j}$ are FIOs with the same
phase $\Phi$ and symbols $\sigma_{k,j}$
belonging to bounded subset of
$S^{-d/2}_{1,0}$, supported in
\begin{equation}\label{b1}
\{(x,\eta)\in\Omega'\times\Gamma:\
|\nabla_x\Phi(x,\eta)|\leq2,
2^{j-1}\leq|\eta|\leq
2^{j+1}\},\quad {\rm if}\
k=0,
\end{equation}
and in
\begin{equation}\label{b2}
\{(x,\eta)\in\Omega'\times\Gamma:\
2^{k-1}\leq|\nabla_x\Phi(x,\eta)|\leq2^{k+1},
2^{j-1}\leq|\eta|\leq
2^{j+1}\},\quad {\rm if}\
k\geq1.
\end{equation}
The operators  $R_{k,j}$ are
smoothing operators whose
symbols $r_{k,j}$ are
 in a bounded subset of
 $\mathcal{S}(\R^{2d})$ and
 supported where
$2^{j-1}\leq|\eta|\leq 2^{j+1}$.\par
Observe that, by the Euler's identity and \eqref{i1},
\[
|\nabla_x\Phi(x,\eta)|=
|\partial^2_{x,\eta}\Phi(x,\eta)\eta|\asymp|\eta|,\quad
\forall (x,\eta)\in\Omega\times\Gamma.
\]
Inserting this equivalence in
\eqref{b1} and \eqref{b2}, we
obtain that there exists
$N_0>0$ such that
$\sigma_{k,j}$ vanishes
identically if $|j-k|>N_0$.\\
On the other hand it is
easily seen that the
 functions
$r_k(x,\eta):=\sum_{j=1}^\infty
r_{j,k}(x,\eta)$ are well
defined and belong to a
bounded set of the Schwartz
space $\mathcal{S}(\R^d)$. As
a consequence, we have
\[
\sum_{j=1}^\infty\psi_k(D)T_{j}=\sum_{j\geq1:
|j-k|\leq N_0}
S_{k,j}+r_k(x,D).
\]
It is then clear that
 the right-hand side
in \eqref{continuazione} can
be dominated by the $\Fur
L^\infty$ norm of $f$ by
using \eqref{psitsr},
\eqref{stime} with
$p=\infty$, and the fact that
$\psi_k(D)$, $R_{j,k}$ and
$r_k(x,D)$ are uniformly
bounded on $\Fur
L^\infty$.\par This concludes
the proof of Theorem
\ref{maintheorem2}.

\section*{Appendix}
In this appendix we give some details
of the proof of \eqref{daprovare}. We
use the same notation as in the proof
of Theorem \ref{maintheorem2}. We also
observe that all the formulae below are
meant to hold on the support of
$\chi^{\nu}_j(u(y,\eta))\sigma_j(\eta,y)$.
\par
First we notice that the
following key formula holds:
\begin{equation}\label{app1}
\nabla_{\eta''}=A(u,v)\nabla_v+O(2^{-j/2})\nabla_u.
\end{equation}
Indeed,
\[
\nabla_{\eta''}=A(u,v)\nabla_v+B(u,v)\nabla_u,
\]
for suitable smooth matrix
$A$ and $B$. On the other
hand, because of our choice
of the splitting
$\eta',\eta''$, we have
$B(u^\nu_j,v)=0$ for every
$v$. A Taylor expansion of
$B(u,v)$ with respect to $u$,
around $u=u^\nu_j$, then
shows \eqref{app1}, because
$|u-u^\nu_j|\leq C_1
2^{-j/2}$.\\
Similarly it turns out
\[
\nabla_v=A'(u,v)\nabla_{\eta''}+O(2^{-j/2})\nabla_{\eta'},
\]
which implies the estimate
\begin{equation}\label{app2}
\frac{\partial
\eta'}{\partial
v}=O(2^{-j/2})
\end{equation}
for the Jacobian matrix
$\frac{\partial
\eta'}{\partial v}$.
Moreover, using \eqref{app2}
and a Taylor expansion of
$\eta'$, as a function of
$(u,v)$, around
$(u^{\nu}_j,0)$, one deduces
\begin{equation}\label{app3}
(\eta-\eta^\nu_j)'=O(2^{-j/2}).
\end{equation}
Now, we have to show that
the repeated application of
the first order operators
\begin{equation}\label{operatori}
2^{-j/2}\partial_{\eta'_k},\quad
\partial_{\eta''_l},\qquad
k=1,...,r,\ \ l=1,...,d-r,
\end{equation} on
\[
e^{2\pi
i(\Phi(\eta^\nu_j,y)-\langle\nabla_1\Phi(\eta^\nu_j,y),\eta^\nu_j\rangle+R^{\nu}_j(\eta,y))}
\chi^\nu_j(u(y,\eta))\sigma_j(\eta,y)
\]
yields expressions dominated by $C'
2^{-jr/2}$.\par
 This is of
course true when these
operators fall on
$\sigma_j(\eta,y)$, since
$\sigma$ is a symbol of order
$-r/2$.
\par
When the operators
$2^{-j/2}\partial_{\eta'_k}$
fall on
$\chi^\nu_j(u(y,\eta))$ one
obtains acceptable terms
because the factor $2^{-j/2}$
in front of the derivative
balances the loss in
\eqref{perdita}. The same
happens for the derivatives
$\partial_{\eta''_l}$,
because of \eqref{app1} and
the fact that
$\chi^\nu_j(u(y,\eta))$ is
constant on the (pieces of)
planes $u=const$.\par Hence
it remains to prove that the
repeated application of the
operators in
\eqref{operatori} on
$R^{\nu}_j(\eta,y)$ yields
uniformly bounded
expressions. To this end,
denote  by $P_\xi$ the
orthogonal projection on the
vector space parallel to the
plane $u=const$ which
contains
$\xi:=\eta^\nu_j+t(\eta-\eta^\nu_j)$,
$0\leq t\leq1$. We observe
that
\[
(\eta-\eta^\nu_j)''-P_\xi
\left((\eta-\eta^\nu_j)''\right)=0\quad
{\rm if}\
(\eta-\eta^\nu_j)'=0.
\]
Hence
\begin{equation}
(\eta-\eta^\nu_j)''=E^\nu_j(t,y,\eta)(\eta-\eta^\nu_j)'+
P_\xi
\left((\eta-\eta^\nu_j)''\right),
\end{equation}
for a suitable matrix
$E^\nu_j(t,y,\eta)$ whose
entries are positively
homogeneous of degree $0$
with respect to $y\in\Gamma$,
and are uniformly bounded,
together with their
derivatives. Now, in
\eqref{restotaylor} we
substitute $\eta-\eta^\nu_j$
by
\[
\eta-\eta^\nu_j=(\eta-\eta^\nu_j)'+
(\eta-\eta^\nu_j)''=(\eta-\eta^\nu_j)'+E^\nu_j(t,y,\eta)
(\eta-\eta^\nu_j)'+ P_\xi
\left((\eta-\eta^\nu_j)''\right).
\]
Since, by assumption, the
gradient
$\nabla_1\Phi(\cdot,y)$ is
constant on the planes
$u=const$, we have
\[
P_\xi
\left((\eta-\eta^\nu_j)''\right)\in{\rm
Ker}\left(d^2_1\Phi(\xi,y)\right).
\]
Using the bi-linearity of the
Hessian one sees that
$R^\nu_j(\eta,y)$ can
therefore be written in the
form
\begin{equation}\label{acs}
\langle G^\nu_j(\eta,y)
(\eta-\eta^\nu_j)',(\eta-\eta^\nu_j)'\rangle,
\end{equation}
where the matrix
$G^\nu_j(\eta,y)$ has entries
which are uniformly bounded,
together with their
derivatives, by $C 2^j$
(because $\Phi(\eta,y)$ is
positively homogeneous
of degree 1 in $y$ and here $|y|\asymp 2^j$).\\
In view of \eqref{app3} we
see that the repeated action
of the operators
\eqref{operatori} on the
expression in \eqref{acs}
yields uniformly bounded
terms.\par This concludes the
proof.

\section*{Acknowledgements}
The author would like to
thank Elena Cordero, Luigi
Rodino and Michael Ruzhansky
for fruitful conversations
and comments.

\end{document}